\documentclass[12pt]{article}  

\usepackage{amsmath}
\usepackage{amssymb}
\usepackage{cite}
\usepackage[usenames]{color}
\usepackage{graphicx}          
\usepackage[dvips]{epsfig}    
\newtheorem{definition}{Definition}
\newtheorem{theorem}{Theorem}

\newtheorem{remark}{Remark}

\newtheorem{lemma}{Lemma}
\newtheorem{assumption}{Assumption}

\title{\LARGE \bf Obstacle Avoidance Problem for Second Degree Nonholonomic Systems}

\author{Victoria Grushkovskaya$^{1,3}$ and Alexander Zuyev$^{2,3}$
\thanks{Corresponding author.\newline
$^{1}$Institute of Mathematics, University of W\"{u}rzburg,       97074 W\"{u}rzburg, Germany
        {\tt\small viktoriia.grushkovska@mathematik.uni-wuerzburg.de}.
        \newline
$^{2}$Max Planck Institute for Dynamics of Complex Technical Systems, 39106 Magdeburg, Germany
{\tt\small zuyev@mpi-magdeburg.mpg.de}.\newline
$^{3}$Institute of Applied Mathematics and Mechanics, National Academy of Sciences of Ukraine, 84100 Sloviansk, Ukraine.\newline
This work was  supported in part by the German Research Foundation (GR 5293/1-1) and the State Fund for Fundamental Research of Ukraine (F75/27190).
}}
\date{}

\begin{document}

\maketitle
\thispagestyle{empty}

\begin{abstract}
In this paper, we propose a new control design scheme for solving the obstacle avoidance problem for nonlinear driftless control-affine systems. The class of systems under consideration satisfies controllability conditions with iterated Lie brackets up to the second order. The time-varying control strategy is defined explicitly in terms of the gradient of a potential function. It is shown that the limit behavior of the closed-loop system is characterized by the set of critical points of the potential function. The proposed control design method can be used under rather general assumptions on potential functions, and particular applications with navigation functions are illustrated by numerical examples.
\end{abstract}

\newpage
\section{Introduction}\label{sec_1}
The development of control algorithms for stabilization and motion
planning problems is one of the most important issues in mathematical control theory
which attracts significant theoretical interest and is highly demanded in various
engineering applications. Motion planning problems were intensively studied by many
researchers, and there exist several classical approaches for their solution. Comprehensive reviews of the main approaches are given, e.g., in~\cite{Kolman95,Sic08,Jean14}. However,
the presence of obstacles in the state space significantly increases the complexity of
motion planning.
 Up to our knowledge, the obstacle avoidance problem was studied only for systems of a special structure (cf.~\cite{Kod90,Rim92,Tan01,Rou08,Va08,Va12,Popa,Pa16a}), {or without ensuring asymptotic stability of the target point~(see, e.g.~\cite{Wal15} and references therein)}, while solving these problems for general classes of nonlinear systems {together with guaranteed stability properties} remains a challenging issue.

In this paper, we present a rigorous theoretical analysis of the obstacle avoidance problem and a proof of
asymptotic stability for a rather general class of driftless control-affine systems. Note that a general class of nonholonomic systems was also considered in the recent paper~\cite{Ura17}, where a smooth time-invariant controller was constructed based on the gradient of a  potential function. However, since nonholonomic systems are not stabilizable by continuous state feedback laws~\cite{Bro83},
the above result ensures only the stability in the sense of Lyapunov (but not asymptotic stability).  This paper continues our study of obstacle avoidance problems for nonholonomic systems. In the previous work~\cite{ZuGrIFAC}, we have proposed time-varying feedback controls ensuring a collision-free motion
for systems of the degree of nonholonomy 1. The contribution of the present paper is twofold.
On the one hand, we present a constructive solution to the obstacle avoidance problem for nonholonomic systems under the controllability condition with second-order iterated Lie brackets. We construct a control algorithm which ensures the collision-free motion of this class of systems. On the other hand, a novel control design scheme is derived for the stabilization of nonholonomic systems. The proposed class of controls is given by
trigonometric polynomials with state-dependent coefficients. Although the exploited approach is similar to the one already presented in~\cite{ZuGrIFAC}, our new formulas for the control coefficients are much simpler {and do not require solving a system of cubic  equations}. This is of particular use for both the motion planning with obstacle avoidance and for further development of control algorithms for nonholonomic systems of an arbitrary high degree of nonholonomy.

The rest of this paper is organized as follows. In the remaining part of this section, we formulate the obstacle avoidance problem and introduce relevant notations and definitions. The main result is  stated in Section~II and illustrated with several examples in Section~III. The proof of the main result is given in the Appendix.

\subsection{Problem statement \& Main idea}
Consider a control system
\begin{equation}
\dot x = \sum_{i=1}^m u_i f_i (x), \quad x\in D\subset \mathbb R^n,\;u\in \mathbb R^m
\label{Sigma}
\end{equation}
where $x=(x_1,...,x_n)^T$ is the state and $u=(u_1,...,u_m)^T$ is the control, $ m<n$, and $f_i\in C^1(D)$.

{To formulate the obstacle avoidance problem, we assume that there is a closed bounded domain $\mathcal W\subset \mathbb R^n$ (workspace) and open domains 
$$\mathcal O_1,\mathcal O_2,..., \mathcal O_N\subset \mathcal W\;\;\text{(obstacles)}$$  such that
$
D=\mathcal W\setminus\bigcup_{j=1}^N\mathcal O_j.
$
We refer to $D$ as the free space and assume that  $D$ is \emph{valid}~\cite{Kod90}, i.e.
$\displaystyle
\overline{\mathcal O_i} \subset {\rm int}\, \mathcal W,\;
\overline{\mathcal O_i} \cap \overline{\mathcal O_j} = \emptyset\;\;\text{if}\; \neq j$,
for all $i,j\in\{1,\dots,N\}$.
Here and in the sequel, we denote the interior of $\mathcal W$ as ${\rm int}\, \mathcal W$ and
its closure as $\bar {\mathcal W}$.

We will study the following obstacle avoidance problem: \emph{for a given initial point $x^0\in {\rm int}\, D$ and a destination point $x^*\in {\rm int}\, D$, our goal is to construct an admissible control such that the corresponding solution of~\eqref{Sigma}
with the initial data $x(0)=x^0$ satisfies the conditions:}
$$
1)\  x(t)\in {\rm int}\, D\text{ for all }t\ge 0;\quad
2)\ x(t)\to x^* \text{ as } t\to +\infty.
$$

\emph{The main objective} of this paper is to solve  the obstacle avoidance problem for system~\eqref{Sigma} by approximating a gradient-like dynamics corresponding to certain potential function. In particular, we aim to construct control laws for~\eqref{Sigma} such that the trajectories of~\eqref{Sigma}  approximate the solution $\bar x(t)$ of
\begin{equation}
 \dot {\bar x} = - \nabla P(\bar x),
\label{grad_system}
\end{equation}
where $\nabla P(\bar x)=\big(\tfrac{\partial P(\bar x)}{\partial \bar x_1},\dots,\tfrac{P(\bar x)}{\partial\bar x_n}\big)^T$  is the gradient of a potential function ensuring the collision free motion of system~\eqref{grad_system}. In particular,  proper candidates for $P$ are navigation functions~\cite{Kod90,Rim92}, artificial potential fields~\cite{Kha86}, etc.

{It should be emphasized that, as the control system~\eqref{Sigma} is underactuated ($m<n$), an arbitrary solution $\bar x(t)$ of~\eqref{grad_system} with the initial condition $\bar x(0)=x^0$ does not satisfy~\eqref{Sigma} in general.
However, the curve $\Gamma=\{\bar x(t):\, 0\le t < \infty\}$ can be approximated by admissible trajectories of~\eqref{Sigma} {with}  high-frequency high-amplitude (open-loop) controls
(see, e.g.,~\cite{Liu97,Jean14,ZuGrBeECC} and references therein).} In this paper, we will present \emph{a novel explicit control design scheme} for solving the obstacle avoidance problem by using time-varying feedback controls and exploiting the sampling concept.

\subsection{Preliminaries}

\paragraph{Notations and definitions}
for a function $f:\mathbb R^n\to\mathbb R$ and a constant $c\in\mathbb R$, we denote $\mathcal L_c=\{x\in\mathbb R^n:f(x)\le c\}$.
 For  vector fields $f,g:\mathbb R^n\to\mathbb R^n $ and $x\in\mathbb R^n$, we denote the Lie derivative as
 $L_gf(x)=\lim\limits_{s\to0}\tfrac{f(x+sg(x))-f(x)}{s}$, and  $[f,g](x)=L_fg(x)-L_gf(x)$ is the Lie bracket. { For any sets $D_1,D_2\subset\mathbb R^n$,  $\rho(D_1,D_2)=\inf_{x\in D_1,y\in D_2}\|x-y\|$.}
Similar to the approaches of~\cite{Clar97,ZuSIAM}, we will exploit the sampling concept.
For a given $\varepsilon>0$, define a partition $\pi_\varepsilon$ of $[0,+\infty)$ into intervals
$
[t_j,t_{j+1})$, $t_j=\varepsilon j$, $j=0,1,2,\dots .$
\begin{definition}
\emph{Given a feedback $u=h(t,x)$, $h:[0,+\infty)\times D\to\mathbb R^m$, $\varepsilon>0$, and $x^0\in\mathbb R^n$, a $\pi_\varepsilon$-solution of~\eqref{Sigma} corresponding to $x^0\in D$ and $h(t,x)$ is an absolutely continuous function  $x(t)\in D$, defined for $t\in[0,+\infty)$, such that  $x(0)=x^0$ and
$
\dot x(t)=f\big(x(t), h(t,x(t_j))\big)$, $t\in[t_j,t_{j+1}),
$
for each j=0,1,2,\dots.}
\end{definition}

\section{Main result}
We focus on the class of systems~\eqref{Sigma} whose vector fields together with their iterated Lie brackets up to the second order span the whole $\mathbb R^n$. Namely, we assume that
 there exist sets of indices $S_1\subseteq \{1,2,...,m\}$,  $S_2\subseteq \{1,2,...,m\}^2$, $S_3 \subseteq\{1,2,...,m\}^3$ such that $|S_1|+|S_2|+|S_3|=n$, and
\begin{equation}\label{rank}
\begin{aligned}
 {\rm span}\big\{f_{i}(x), &[f_{j_1},f_{j_2}](x), {\left[[f_{\ell_1},f_{\ell_2}],f_{\ell_3}\right](x)} \,|\,\\
 & i{\in}S_1,
(j_1,j_2){\in} S_2, (\ell_1,\ell_2,\ell_3){\in} S_3\big\}{=}\mathbb{R}^n,
\end{aligned}
\end{equation}
for each $x\in D$.
 Under the above assumption, the following $n\times n$-matrix is nonsingular for each $x\in D $:
\begin{equation}\label{fmatrix}
\begin{aligned}
 \mathcal F(x)= \Big(\big(f_{i}(x)&\big)_{j_1\in S_1}\ \big([f_{j_1},f_{j_2}](x)\big)_{(j_1,j_2)\in S_2} \\
 &\big(\left[[f_{\ell_1},f_{\ell_2}],f_{\ell_3}\right](x)\big)_{(\ell_1,\ell_2,\ell_3)\in S_3}\Big).
\end{aligned}
\end{equation}
  Below we introduce a  family of control functions which will be used for
  approximating the gradient dynamics corresponding to a function $P\in C^1(D;\mathbb R)$ by the trajectories of~\eqref{Sigma}.
  Namely, for given positive real number $\varepsilon$ and $\nu$, we take
\begin{align}
 u&^\varepsilon_k(t,x)=\sum_{i_1\in S_1} a_{i_1}(x)\phi^{(k,\varepsilon)}_{i_1}(t)\nonumber \\
 &+\varepsilon^{-\tfrac{1}{2}}\sum_{(j_1,j_2)\in S_2} \sqrt{|a_{j_1j_2}(x)|}\phi^{(k,\varepsilon)}_{j_1j_2}(t)\label{cont}\\
  &+\varepsilon^{-\tfrac{2}{3}}\sum_{(\ell_1,\ell_2,\ell_3)\in S_3} \sqrt[3]{a_{\ell_1\ell_2\ell_3}(x)}\phi^{(k,\varepsilon)}_{\ell_1\ell_2\ell_3}(t),k=1,\dots,m,\nonumber
\end{align}
and define the  state dependent vector function
\begin{align}
a(x)&=\big(a_{i_1}(x)\big|_{i_1\in S_1}, a_{j_1j_2}(x)\big|_{(j_1,j_2)\in S_2}, \nonumber\\
&\qquad\qquad\qquad a_{\ell_1\ell_2\ell_3}(x)\big|_{(\ell_1,\ell_2,\ell_3)\in S_3}\big)^T\in\mathbb R^n\nonumber\\
&=- \gamma \mathcal F^{-1}(x) \nabla P(x)\label{a}
\end{align}
with some $\gamma>0$,  and
\begin{align}
\phi^{(k,\varepsilon)}_{i_1}(t)&{=}\delta_{k i_1},\nonumber \\
 \phi^{(k,\varepsilon)}_{j_1j_2}(t)&{=}2\sqrt{\pi K_{j_1j_2}}\Big(\delta_{kj_1}{\rm sign}(a_{j_1,j_2}(x))\cos{\frac{2\pi K_{j_1j_2}}{\varepsilon}}t\nonumber\\
 &{+}\delta_{kj_2}\sin{\frac{2\pi K_{j_1j_2}}{\varepsilon}}t\Big), \nonumber\\
\phi^{(k,\varepsilon)}_{\ell_1\ell_2\ell_3}(t)&{=}2\sqrt[3]{2\pi^2 K_{3\ell_1\ell_2\ell_3}K_{4\ell_1\ell_2\ell_3}}\Big(\delta_{kl3_1}\cos{\frac{2\pi K_{1\ell_1\ell_2\ell_3}t}{\varepsilon}}\nonumber\\
&{+}\delta_{k\ell_2}\sin{\frac{2\pi K_{2\ell_1\ell_2\ell_3}t}{\varepsilon}} \nonumber\\
&{+}\delta_{k\ell_3}{\cos{\frac{2\pi K_{1\ell_1\ell_2\ell_3}t}{\varepsilon}}\sin{\frac{2\pi K_{2\ell_1\ell_2\ell_3}t}{\varepsilon}}}\Big).\label{phi}
\end{align}
Here $\delta_{ki}$ is the Kronecker delta, and non-zero integer parameters
$K_{j_1j_2}$, $K_{1\ell_1\ell_2\ell_3}$, $K_{2\ell_1\ell_2\ell_3}$, $ K_{3\ell_1\ell_2\ell_3}=K_{1\ell_1\ell_2\ell_3}+K_{2\ell_1\ell_2\ell_3}$, and $K_{4\ell_1\ell_2\ell_3}=K_{2\ell_1\ell_2\ell_3}-K_{1\ell_1\ell_2\ell_3}$ are specified according to the following non-resonance assumption.
\begin{assumption}~\label{ControlDesign}
\emph{
\begin{itemize}
  \item $|K_{i_1i_2}|\ne |K_{j_1j_2}|$ for all $S_2\ni (i_1,i_2)\ne(j_1,j_2)\in S_2$;
  \item  for any $s_1,s_2,s_3{\in}\{1,\dots,4\}$, if  $c_1,c_2,c_3{\in}\{+1,-1\}$ and
  $$c_1K_{s_1i_1i_2i_3}{+}c_2K_{s_2j_1j_2j_3}{=}c_3K_{s_3\ell_1\ell_2\ell_3},$$
  then $s_1\ne s_2{\ne} s_3{\ne} s_1$ and $(i_1,i_2,i_3){=}(j_1,j_2,j_3){=}(\ell_1,\ell_2,\ell_3){\in} S_3$.
\end{itemize}}
  \end{assumption}
The aim of the above non-resonance assumption is to ensure the motion of system~\eqref{Sigma} in the direction of the Lie brackets appearing in rank condition~\eqref{rank} without generating undesirable equivalent Lie brackets, as will be discussed later.

The following theorem shows that, under a proper choice of the potential function $P$, the control functions defined in~\eqref{cont} solve the obstacle avoidance problem for system~\eqref{Sigma}.
\begin{theorem}~\label{thm_main}
 \emph{Let $D\subset \mathbb R^n$ be a bounded free space, $f_i\in C^3(D;\mathbb R)$, $i=1,\dots,m$, and let there exist an  $\alpha>0$ such that
 $
    \|\mathcal F^{-1}(x)\|\le \alpha \text{ for all }x\in D,
$
where the matrix $\mathcal F(x)$ is given by~\eqref{fmatrix}. Suppose also that  a function $P\in C^2(D)$ is such that $\rho(D_0,\partial D)>0$, for any $x^0\in D$ and $D_0:=\mathcal L_{P(x^0)}=\{x\in\mathbb R^n:P(x)\le P(x^0)\}$, and let the functions $u^\varepsilon_k(t,x)$, $k=1,\dots,m$, be defined in~\eqref{cont}--\eqref{phi} with the parameters satisfying Assumption~\ref{ControlDesign}.}

\emph{Then  there exists an $\bar\varepsilon>0$  such that, for any $\varepsilon\in(0,\bar\varepsilon]$, the $\pi_\varepsilon$-solution of system~\eqref{Sigma} with the controls $u_k=u^\varepsilon_k(t,x)$ and the initial data $x(0)=x^0\in D$ is well-defined on $t\in [0,+\infty)$ and satisfies the following properties:
$$
 x(t) \in {\rm int}\, D,\quad \forall t\ge 0,
 $$
 $$
 x(t) \to Z_0 = \{x\in D_0\,:\,\nabla P (x)=0\}\quad \text{as}\;\; t\to +\infty.
$$}
\end{theorem}
The proof of Theorem~1 is in the Appendix.

\begin{remark}
  It is important that the matrix $\mathcal F(x^0)$ contains only the Lie brackets of $f_i$ which appear in the controllability condition~\eqref{rank}.
    For example, the motion in the direction of the Lie bracket $[[f_1,f_2],f_3](x)$ can be generated with the controls
  $$
 \begin{aligned}
  &u_1^\varepsilon(t)=\cos\frac{2\pi K_{1123}t}{\varepsilon},\, u_2^\varepsilon(t)=\sin\frac{2\pi K_{2123}t}{\varepsilon},\\
  &u_3^\varepsilon(t)=\cos\frac{2\pi K_{1123}t}{\varepsilon}\sin\frac{2\pi K_{2123}t}{\varepsilon}.
 \end{aligned}
 $$
  In this case, the solutions of system~\eqref{Sigma} can be represented as
 follows: $$x(\varepsilon)=\frac{\varepsilon^3}{16\pi^2(K_{2123}^2-K_{1123}^2)}{\left[\big[f_{1},f_{2}],f_{3}\right](x^0)}+R(\varepsilon),$$ provided that $|K_{1123}|\ne|K_{2123}|$, $|K_{1123}|\ne2|K_{2123}|$, and $2|K_{1123}|\ne |K_{2123}|$. {Here $R(\varepsilon) $ denotes higher order terms, $ \|R(\varepsilon)\|\le \sigma \varepsilon^{4/3}\|\nabla P(x^0)\|^{4/3}$ with some $\sigma>0$ (see the proof of Theorem~\ref{thm_main}) .} Similarly, the controls
 $$
\begin{aligned}
 & u_1^\varepsilon(t)=\cos\frac{2\pi K_{1123}t}{\varepsilon}\left(1+\sin\frac{2\pi K_{2123}t}{\varepsilon}\right),\\
 & u_2^\varepsilon(t)=\sin\frac{2\pi K_{2123}t}{\varepsilon},\,u_3^\varepsilon(t)=0,
\end{aligned}
  $$
  generate the only Lie bracket $[[f_1,f_2],f_1](x^0)$ (with higher order terms being omitted).
\end{remark}

\section{Example}
In this section, we illustrate the proposed control design approach with several second degree nonholonomic systems.
\paragraph{Rigid body with oscillators}
Consider the system
\begin{equation}\label{bro}
\begin{aligned}
& \dot x_1=u_1,\; \dot x_2=u_2,\;  \dot x_3=x_1^2u_2-x_2^2u_1.
\end{aligned}
\end{equation}
System~\eqref{bro} has been considered in~\cite{YKD} as a mathematical model of a rotating rigid body with oscillating point masses.
It is easy to check that the above system satisfies the rank condition~\eqref{rank} with $S_1=\{1,2\}$, $S_2=\emptyset$, $S_3=\{(1,2,1)\}$, i.e.
$
{\rm rank}\ \mathcal F(x)={\rm rank}\ \left(f_1(x)\ f_2(x)\ \big[[f_1,f_2],f_1\big]{(x)}\right)=3
$
for each $x\in\mathbb R^3$.
According to Theorem~1, we define
\begin{equation}\label{cont_bro}
\begin{aligned}
u_1^\varepsilon(t,x)=&a_1(x)+2\sqrt[3]{\frac{2\pi^2(K_{2121}^2-K_{1121}^2)}{\varepsilon^2}}\sqrt[3]{a_{121}(x)}\\
&\times\cos\tfrac{2\pi K_{1121}t}{\varepsilon}\Big(1+\sin\tfrac{2\pi K_{2121}t}{\varepsilon}\Big),\\
u_2^\varepsilon(t,x)=&a_2(x)+2\sqrt[3]{\frac{2\pi^2(K_{2211}^2-K_{1121}^2)}{\varepsilon^2}}\sqrt[3]{a_{121}(x)}\\
&\times\sin\tfrac{2\pi K_{2121}t}{\varepsilon},
\end{aligned}
\end{equation}
where
$$
\begin{aligned}
\left(
\begin{aligned}
 &a_1(x)\\
 &a_2(x)\\
 &a_{12}(x)
\end{aligned}
\right)
&=
-\gamma\mathcal F^{-1}(x)\nabla P(x)\\
&=\left(
                                       \begin{aligned}
                                         &-\gamma \tfrac{\partial P(x)}{\partial x_1}\\
                                         &-\gamma \tfrac{\partial P(x)}{\partial x_2}\\
                                         &\frac{\gamma}{2}\big(x_2^2\tfrac{\partial P(x)}{\partial x_1}-x_1^2\tfrac{\partial P(x)}{\partial x_2}-\tfrac{\partial P(x)}{\partial x_3}\big)
                                       \end{aligned}
                                     \right).
\end{aligned}
$$

\begin{figure}[ht]
\begin{center}
  \includegraphics[width=0.5\linewidth]{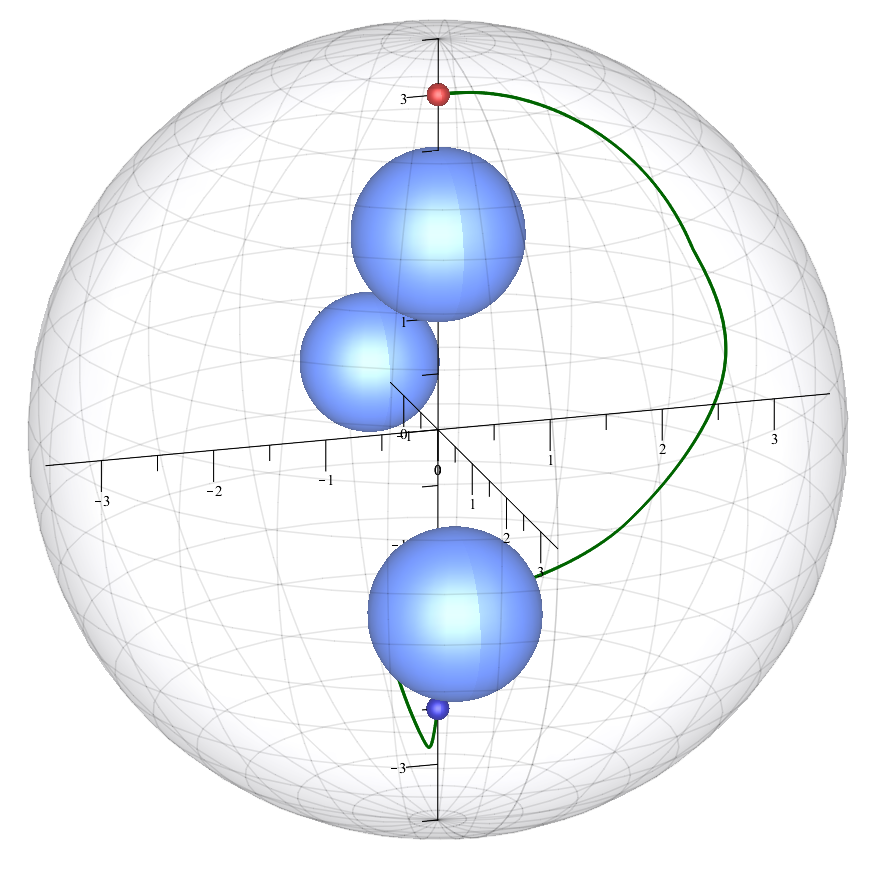}\\
\includegraphics[width=0.5\linewidth]{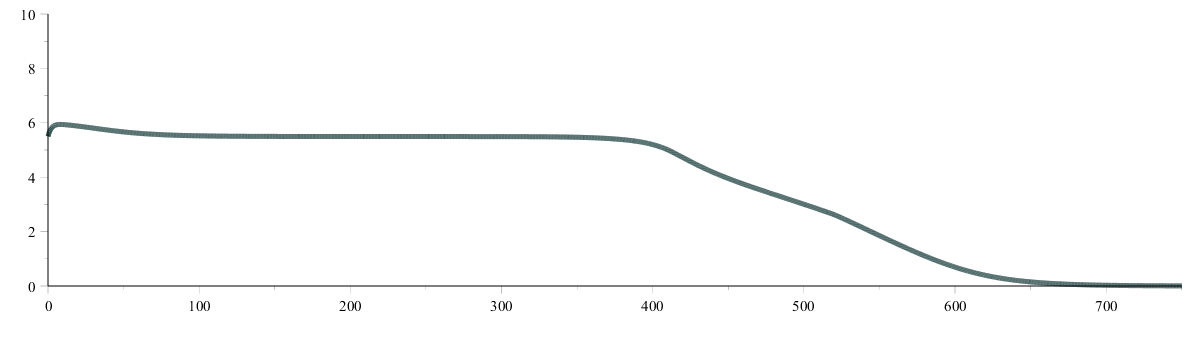}\\
\includegraphics[width=0.5\linewidth]{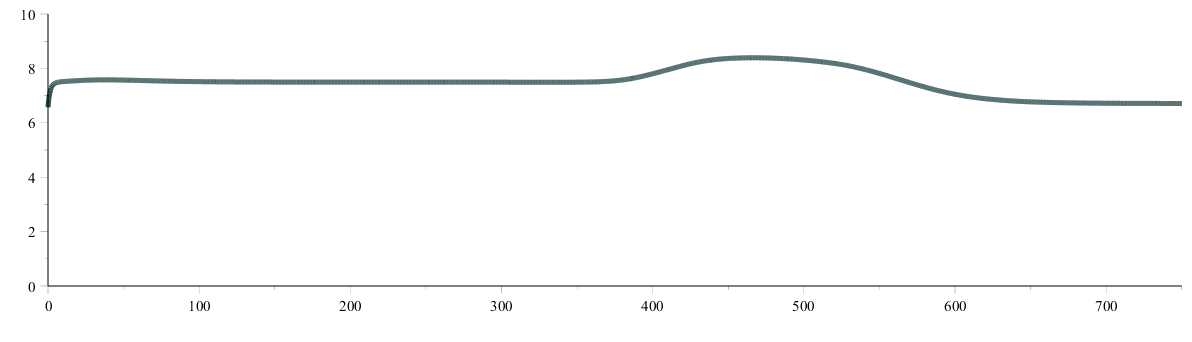}
\end{center}
\caption{Trajectory of system~\eqref{bro} with controls~\eqref{cont_bro} (top), time-plot of the function $\|x(t)-x^*\|$ (middle), and $\ln(1+\prod\limits_{j=0}^3\beta_j(x(t)))$ (bottom). }
\end{figure}
In this example,  we put $K_{1121}=1$, $K_{2121}=3$, $\varepsilon=0.5$, $\gamma=0.5$. It is easy to see that Assumption~1 is satisfied.
Suppose that the goal is to steer system~\eqref{bro} to  the target $x^*=(0,0,3)^T$  avoiding collisions with three spherical obstacles ${\mathcal O}_j=\{x\in \mathbb R^3: \beta_j(x)<0\}$ and remaining in the workspace $\mathcal W = \{x\in \mathbb R^3: \beta(x)\ge 0\}$, where$$
\begin{aligned}
&\beta_0(x)=12.25-x_1^2-x_2^2-x_3^2,\\
&\beta_1(x)=x_1^2+x_2^2+(x_3-1.75)^2-0.5625,\\
&\beta_2(x)=(x_1-0.5)^2+x_2^2+(x_3+1.5)^2-0.5625,\\
&\beta_3(x)=(x_1+2)^2+x_2^2+x_3^2-0.36.
\end{aligned}
$$
 As a potential function $P(x)$, we take the following navigation function proposed in~\cite{Kod90}:
\begin{equation}\label{p_nav}
P(x)=\frac{\|x-x^*\|^2}{\big(\|x-x^*\|^{4}+\prod\limits_{j=0}^3\beta_j(x)\big)^{1/2}}.
\end{equation}
The behavior of the corresponding closed-loop system with the initial point $x^0=(0,0,-3)^T$ is illustrated in Fig.~1. The time-plot of the function $\ln(1+\prod\limits_{j=0}^3\beta_j(x(t)))$ illustrates that the trajectory of system~\eqref{bro} remains in the free space.

\paragraph{Disc rolling on the plane}
Consider the control system that describes the motion of a unit disc rolling on the horizontal plane~(cf. \cite{LC90}):
\begin{equation}\label{disc}
 \dot x_1=u_1\cos x_3,\; \dot x_2=u_1\sin x_3,\; \dot x_3=u_2,\;  \dot x_4=u_1.
\end{equation}
The vector fields of~\eqref{disc} satisfy~\eqref{rank} with $S_1=\{1,2\}$, $S_2=\{(1,2)\}$, $S_3=\{(1,2,2)\}$. In this case, we take
\begin{equation}\label{cont_disc}
\begin{aligned}
&u_1^\varepsilon(t,x)=a_1(x)+2\sqrt{\frac{\pi K_{12}}{\varepsilon}\cos\tfrac{2\pi K_{12}t}{\varepsilon}}\\
&+2\sqrt[3]{\frac{2\pi^2(K_{2122}^2-K_{1122}^2)}{\varepsilon^2}}\sqrt[3]{a_{122}(x)}\cos\tfrac{2\pi K_{1122}t}{\varepsilon},\\
&u_2^\varepsilon(t,x)=a_2(x)+2\sqrt{\frac{\pi K_{12}}{\varepsilon}\sin\tfrac{2\pi K_{12}t}{\varepsilon}}\\
&+2\sqrt[3]{\frac{2\pi^2(K_{2122}^2-K_{1122}^2)}{\varepsilon^2}}\sqrt[3]{a_{122}(x)}\\
&\times\sin\tfrac{2\pi K_{2122}t}{\varepsilon}\Big(1+\cos\tfrac{2\pi K_{1122}t}{\varepsilon}\Big),
\end{aligned}
\end{equation}
with
$$
\begin{aligned}
\left(
\begin{aligned}
 &a_1(x)\\
 &a_2(x)\\
 &a_{12}(x)\\
&a_{122}(x)
\end{aligned}
\right)
{=}\left(
                                       \begin{aligned}
                                         &-\gamma \tfrac{\partial P(x)}{\partial x_4}\\
                                         &-\gamma \tfrac{\partial P(x)}{\partial x_3}\\
                                         &-{\gamma}\big(\sin x_3\tfrac{\partial P(x)}{\partial x_1}+\cos x_3\tfrac{\partial P(x)}{\partial x_2}\big)\\
&\gamma\big(\cos x_3\tfrac{\partial P(x)}{\partial x_1}+\sin x_3\tfrac{\partial P(x)}{\partial x_2}+\tfrac{\partial P(x)}{\partial x_4}\big)
                                       \end{aligned}
                                     \right).
\end{aligned}
$$

\begin{figure}[ht]
\begin{center}
  \includegraphics[width=0.5\linewidth]{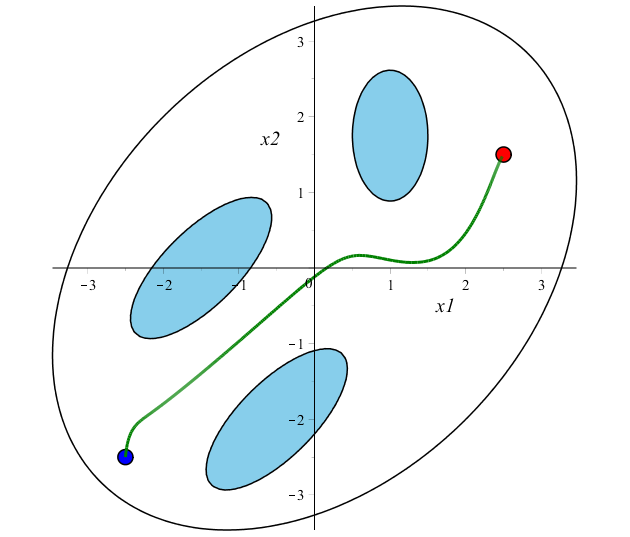}\\
\includegraphics[width=0.5\linewidth]{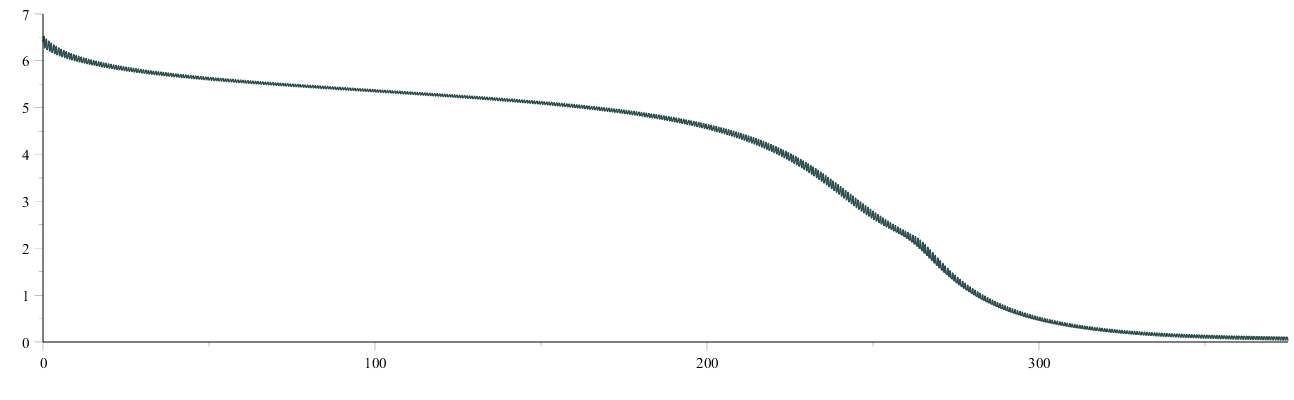}\\
\includegraphics[width=0.5\linewidth]{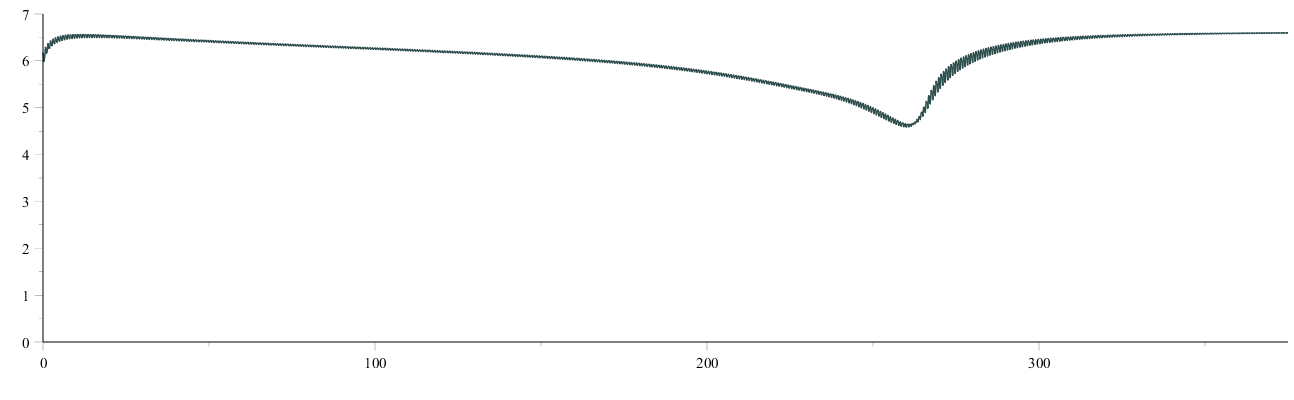}
\end{center}
\caption{Projection of the trajectory of system~\eqref{disc} with controls~\eqref{cont_disc} on the  $(x_1,x_2)$-plane (top), time-plot of the function $\|x(t)-x^*\|$ (middle), and $\ln(1+\prod\limits_{j=0}^3\beta_j(x(t)))$ (bottom). }
\end{figure}
The function $P$ is defined by~\eqref{p_nav}, and the integer parameters $K_{12}=1$, $K_{1122}=3$, $K_{2122}=7$ satisfy Assumption~1.
For the simulation, we take  the initial point  $x^0=\big({-2.5},-2.5,0,\tfrac{\pi}{4}\big)^T$, the target point  $x^*=\big(2.5,1.5,0,\tfrac{\pi}{4}\big)^T$, $\varepsilon=0.75$, $\gamma=0.5$,  $\lambda=2$, while the obstacles and the workspace are defined by
$$
\begin{aligned}
&\beta_0(x)=1-\tfrac{1}{8}(x_1-x_2)^2-\tfrac{1}{32}(x_1+x_2)^2-\tfrac{1}{10}x_3^2-\tfrac{1}{10}x_4^2,\\
&\beta_1(x)=2(x_1-x_2-1.5)^2+\tfrac{1}{3}(x_1+x_2+2.5)-1,\\
&\beta_2(x)=2(x_1-x_2+1.5)^2+\tfrac{1}{3}(x_1+x_2+1.5)-1,\\
&\beta_3(x)=4(x_1-1)^2+\tfrac{4}{3}(x_2-1.75)-1.
\end{aligned}
$$
Fig.~2 illustrates that the controller proposed can be used for solving the obstacle avoidance problem.

\section{Conclusions}

The main theoretical contribution of this paper is summarized in Theorem~1 and illustrated with numerical examples in Section~III. It should be emphasized that the proposed control design scheme is obtained for the general class of second degree nonholonomic systems and  is applicable under rather general assumptions on the potential function $P(x)$.
The advantages of our approach are twofold in nature. On the one hand, the controllers are defined explicitly by  rather simple formulas~\eqref{cont}--\eqref{phi}.
       These formulas exploit the inversion of the matrix composed by the vector fields appearing in the controllability conditions.
 On the other hand, we use a family of oscillating controls $u^\varepsilon(t,x)$ with state-depending coefficients in contrast to the open-loop approaches to motion planning of nonholonomic systems. {We expect that a similar control design scheme with the coefficients defined by~\eqref{a} can be extended to general nonholonomic systems. We leave this issue for future work.}

\section*{APPENDIX}
We will use the following result.
\begin{lemma}~{\rm \bf \cite{ZuSIAM}}~\label{lemma_x}
Let $x(t)\in  D\subset \mathbb R^n$, $0\le t\le \tau$, be a solution of system~\eqref{Sigma} with a control $u\in C[0,\tau]$, and let
$$
\|f_i(x)-f_i(y)\|\le L \|x-y\|,\quad \forall x,y\in  D,\; i= 1,..., m.
$$
with some $L>0$.
Then
\begin{equation}
\|x(t)-x(0)\|\le \frac{M}{L}(e^{LU t}-1),\quad t\in [0,\tau],
\label{apriori_est}
\end{equation}
where $M=\sup\limits_{\hspace{-1em}\underset{1\le i\le m}{x\in D}}\|f_i(x)\|$, $U=\max\limits_{0\le t\le \tau}\sum_{i=1}^m |u_i(t)|$.
\end{lemma}
{\em Proof of Theorem~\ref{thm_main}}. Without loss of generality, we assume that $x^*=0$. The proof is splitted into several steps.

  \emph{Step~1.} Let $x^0\in D$, and let $P(x)$ be a given potential function.  On the first step we  show that there exists an $\varepsilon_1>0$ such that, for all $\varepsilon\in(0,\varepsilon_1)$, the $\pi_\varepsilon$-solution of system~\eqref{Sigma} with the initial data $x(0)=x^0$ and the controls $u_k=u^\varepsilon_k(t,x)$ defined by~\eqref{cont} is well-defined on $t\in [0,\varepsilon]$.
Let us define $D_0{ =} {\mathcal L}_{P (x^0)}\subset D$, $d_0 = \rho\left(D_0,\partial D \right)>0$, $d_1\in (0,d_0)$,
$$
M=\sup\limits_{\hspace{-1em}\underset{1\le i\le m}{x\in D_0}}\|f_i(x)\|,\,U(x^0)=\max\limits_{0\le t\le \varepsilon}\sum_{i=1}^m |u_i^\varepsilon(t,x^0)|.
$$
Using H\"{o}lder's inequality, one can estimate $U$ as
\begin{align}
&U(x^0)\le\sum_{j_1\in S_1}|a_{j_1}(x^0)| +4\sum_{(j_1,j_2)\in S_2}\sqrt{\frac{\pi K_{j_1j_2}|a_{j_1j_2}(x^0)|}{\varepsilon}} \nonumber\\
&+6\sqrt[3]{\frac{2\pi^2}{\varepsilon^2}}\sum_{(\ell_1,\ell_2,\ell_3)\in S_3}\sqrt[3]{|K_{2\ell_1\ell_2\ell_3}^2{-}K_{1\ell_1\ell_2\ell_3}^2||a_{\ell_1\ell_2\ell_3}(x^0)|}\nonumber\\
&\le \sqrt{|S_1|}a_0+4\sqrt{\frac{\pi}{\varepsilon}}\Big(\sum_{(j_1,j_2)\in S_2}{K_{j_1j_2}}^{2/3}\Big)^{3/4}{a_0^{1/2}}\label{est_u}\\
&+6\sqrt[3]{\frac{2\pi^2}{\varepsilon^2}}\Big(\sum_{(\ell_1,\ell_2,\ell_3)\in S_3}|K_{2\ell_1\ell_2\ell_3}^2-K_{1\ell_1\ell_2\ell_3}^2|^{2/5}\Big)^{5/6}{a_0^{1/3}}\nonumber\\
&\le C_1\|\nabla P(x^0)\|{+}\frac{C_2}{\sqrt\varepsilon}\|\nabla P(x^0)\|^{1/2}{+}\frac{C_3}{\sqrt[3]{\varepsilon^2}}\|\nabla P(x^0)\|^{1/3},\nonumber
\end{align}
where $a_0=\|a(x^0)\|$, $C_1= \gamma\alpha\sqrt{|S_1|}$,

 $C_2=4\sqrt{{\pi\gamma\alpha}}\Big(\sum\limits_{(j_1,j_2)\in S_2}{K_{j_1j_2}}^{2/3}\Big)^{3/4}$,

 $C_3{=}6\sqrt[3]{2\pi^2\gamma\alpha}\Big(\sum\limits_{(\ell_1,\ell_2,\ell_3)\in S_3}|K_{2\ell_1\ell_2\ell_3}^2{-}K_{1\ell_1\ell_2\ell_3}^2|^{2/5}\Big)^{5/6}$.

Let $M_P=\sup_{x\in D}\|\nabla P(x)\|$. Note that $M_P$ is well-defined as $D$ is compact and $\nabla P$ is continuous in $D$. Then it follows from Lemma~\ref{lemma_x} and estimate~\eqref{est_u} that
$$
\begin{aligned}
\|x(t)&-x^0\|\le \frac{M}{L}(e^{U(x^0) L\varepsilon}{-}1)\\
&\le  \frac{M}{L}(e^{L(\varepsilon M_PC_1+\sqrt{\varepsilon M_P}C_2+\sqrt[3]{\varepsilon M_P}C_3)}{-}1),\;t{\in}[0,\varepsilon].
\end{aligned}
$$
Thus, defining $\varepsilon_1$ as the smallest positive root of the equation
$$
\varepsilon M_PC_1+\sqrt{\varepsilon M_P}C_2+\sqrt[3]{\varepsilon M_P}C_3=\frac{1}{L}\Big(\ln\frac{Ld_1}{M}+1\Big),
$$
we obtain
$\|x(t)-x^0\|<d \text{ for each }\varepsilon\in(0,\varepsilon_1)\; \text{and}\; t\in[0,\varepsilon],$
that is, $x(t)\in D$ for all  $t\in[0,\varepsilon]$.

 \emph{Step~2.} In this step, we expand the  $\pi_\varepsilon$-solution of system~\eqref{Sigma} into the Volterra series and estimate its remainder. Namely, we compute
\begin{equation}
 \begin{aligned}
&x(t)=x^0+ \sum_{j_1=1}^mf_{j_1}(x^0)\int\limits_0^t u_{j_1}(s_1){\rm d}s_1\\
&+\sum_{j_1,j_2=1}^m L_{f_{j_2}}f_{j_1}(x^0)\int\limits_0^t\int\limits_0^{s_1} u_{j_1}(s_1)u_{j_2}(s_2){\rm d}s_2{\rm d}s_1\\
&+{\sum\limits_{\hspace{-1em}\ell_1,\ell_2,\ell_3=1}^m} L_{f_{\ell_3}} L_{f_{\ell_2}}f_{\ell_1}(x^0){\int\limits_0^t}{\int\limits_0^{s_2}}{\int\limits_0^{s_1}} u_{\ell_1}(s_1)u_{\ell_2}(s_2)\\
&\quad\quad\times u_{\ell_3}(s_3){\rm d}s_3{\rm d}s_2{\rm d}s_1+r(t),
  \end{aligned}~\label{volt}
\end{equation}
 where the remainder of the Volterra series expansion can be represented as follows (see, e.g.,~{\cite{ZuGrIJC}}):
 \begin{equation}\label{rem}
\begin{aligned}
 r(t)=&{\sum\limits_{\hspace{-1em}j_1,\dots,j_4=1}^m}{\int\limits_0^t}\dots{\int\limits_0^{s_3}} L_{f_{j_4}}\dots L_{f_{j_2}}  f_{j_1}(x(s_4))
 \\
 & \times
  u_{j_4}(s_4)\dots u_{j_1}(s_1){\rm d}s_4\dots{\rm d}s_1.
\end{aligned}
\end{equation}
Denote $H=\frac{1}{24}\max\limits_{x\in D_0}{\sum\limits_{\hspace{-.1em}j_1,\dots,j_4=1}^m}\big\|L_{f_{j_4}}\dots L_{f_{j_2}} f_{j_1}(x)\big\|$. Then
\begin{equation}\label{r_est}
     \|r(t)\|\le {H}\big(U(x^0)t\big)^4.
\end{equation}
Computing the integrals in~\eqref{volt} {with regard  to Assumption~1}, we obtain the following representation for $t=\varepsilon$:
 \begin{align*}
x(\varepsilon)&=x^0+\varepsilon\sum_{j_1\in S_1}f_{j_1}(x^0)a_{j_1}(x^0)\\
&+{\varepsilon}\sum_{(j_1,j_2)\in S_2}[f_{j_1},f_{j_2}](x^0){a_{j_1j_2}(x^0)}\\
&+ {\varepsilon}\sum_{(\ell_1,\ell_2,\ell_3)\in S_3}[[f_{\ell_1},f_{\ell_2}],f_{\ell_3}](x^0){a_{\ell_1\ell_2\ell_3}(x^0)}\\
&+\Omega(a,\varepsilon)+r(\varepsilon)={x^0+\varepsilon}\mathcal F(x^0)a(x^0)+\Omega(a,\varepsilon)+r(\varepsilon),
  \end{align*}
where $\Omega(a,\varepsilon)$ contains the higher order terms. We omit the explicit formula
due to lack of space.  It can be shown that there exists a $\varpi\ge0$ such that
$
\|\Omega(a,\varepsilon)\|\le \varepsilon^{4/3}\varpi\|a(x)\|^{4/3}\text{ for all }x\in D.
$
Taking into account~\eqref{est_u} and~\eqref{r_est}, we can estimate $r(\varepsilon)$ as
$$
\begin{aligned}
\|r(\varepsilon)\|\le &{H}\big(U(x^0)\varepsilon\big)^4\le \varepsilon^{4/3}H\|\nabla P(x^0)\|^{4/3}\big(C_1(\varepsilon M_P)^{2/3}\\
&+{C_2}(\varepsilon M_P)^{1/6}+{C_3}\big)^4.
\end{aligned}
$$
Using the formula~\eqref{a} for $a(x)$, we may rewrite the obtained representation for $x(\varepsilon)$
 as
 \begin{equation}\label{x_volt}
 x(\varepsilon)=x^0-\varepsilon\gamma\nabla P(x^0)+R(\varepsilon),
 \end{equation}
 where $R(\varepsilon)=\Omega(a,\varepsilon)+r(\varepsilon)$, and
\begin{equation}\label{R_est}
 \begin{aligned}
 \|R(\varepsilon)\|\le \sigma \varepsilon^{4/3}\|\nabla P(x^0)\|^{4/3}\;\text{ for all }\varepsilon\in(0,\varepsilon_1),
 \end{aligned}
\end{equation}
 with the fixed $\varepsilon_1$ defined in Step~1, and $\sigma=\varpi(\gamma\alpha)^{4/3}+\big(C_1(\varepsilon_1 M_P)^{2/3}
+{C_2}(\varepsilon_1 M_P)^{1/6}+{C_3}\big)^4$.

\emph{Step 3.} The next goal is to ensure that the function $P$ is decreasing along the trajectories of system~\eqref{Sigma}.

For this purpose, we apply Taylor's formula with the Lagrange form of the remainder to  $P(x(\varepsilon))$ and exploit the formula~\eqref{x_volt}:
$$
\begin{aligned}
P&(x(\varepsilon))=P(x^0)+\nabla P(x^0)(x(\varepsilon)-x^0)^T\\
&+\frac{1}{2}\sum_{i,j=1}^m\left.\frac{\partial^2 P(x)}{\partial x_i\partial x_j}\right|_{\theta}(x_i(\varepsilon)-x_i^0)(x_j(\varepsilon)-x_j^0)\\
&\le P(x^0)-\varepsilon\gamma\|\nabla P(x^0)\|^2+\|\nabla P(x^0)\|\|R(\varepsilon)\|\\
&+\mu\big(\varepsilon\gamma\|\nabla P(x^0)\|+\|R(\varepsilon)\|\big)^2,
\end{aligned}
$$
where $\mu=\frac{1}{2}\sup_{x\in D}\sum_{i,j=1}^m\Big\|\frac{\partial^2 P(x)}{\partial x_i\partial x_j}\Big\|$, $\|\theta-x^0\|\le\|x(\varepsilon)-x^0\|$.
With~\eqref{R_est}, we get
$$
\begin{aligned}
P(x(\varepsilon))\le P&(x^0)-\varepsilon\gamma\|\nabla P(x^0)\|^2\big(1-\mu\varepsilon\gamma\\
&-(1+2\mu\varepsilon\gamma)\sigma\varepsilon^{1/3} M_P^{1/3}-\mu\sigma^2\varepsilon^{5/3} M_P^{2/3}\big).
\end{aligned}
$$
Taking $\varepsilon_2$ as the smallest positive root of the equation
$$
\mu\varepsilon\gamma+(1+2\mu\varepsilon\gamma)\sigma\varepsilon^{1/3} M_P^{1/3}+\mu\sigma^2\varepsilon^{5/3} M_P^{2/3}=1,
$$
we obtain that, for all $\varepsilon\in(0,\varepsilon_2)$ and $x^0\in D_0,$
$$
P(x(\varepsilon))<P(x^0),
$$
provided that $\|\nabla P(x^0)\|\ne 0$. Iterating the above procedure for $x^0\in D_0$, we conclude that $x((j+1)\varepsilon)\in D_0$, and
$
x((j+1)\varepsilon)\in D_0 $, $P(x((j+1)\varepsilon))\le P(x(j\varepsilon))$, $j=0,1,...,$  .

\emph{Step 4.} In this final step, we prove the assertion of Theorem~1.
Consider the discrete-time dynamical system
\begin{equation}
x^j = h(x^{j-1}),\quad \; j=1,2,... \;,
\label{discrete_DS}
\end{equation}
where $h:D_0\to D_0$ maps any $ \xi\in D_0$ to the solution of system~\eqref{Sigma} with the initial condition $x|_{t=0}{=}\xi$ and controls~\eqref{cont} evaluated at $t=\varepsilon$, and $h(\xi)=\xi$ if $\nabla P(\xi)=0$. One can see that
 $x^j=x(j\varepsilon)$, $j=0,1,2,...$, where  $x(t)$ is the $\pi_\varepsilon$ solution of system~\eqref{Sigma} with the initial condition $x|_{t=0}{=}\xi$ and controls $u_k^\varepsilon(t,x)$ given by~\eqref{cont}. As it has been already proved in Step~3,   $x(t)$ is well defined on each interval $I_j=[\varepsilon j,\varepsilon (j+1))$, and $x(t)\in {\rm int}\, D$ for all $t\ge 0$.
Then the  invariance principle~\cite{LaSalle}\footnote{Another versions of the invariance principle are presented in~\cite{K09,Zu2001} for problems of partial stability.} for~\eqref{discrete_DS} implies that
\begin{equation}
x^j \to Z_0\quad\text{ as }\;\; j \to +\infty,
 \label{limit_discrete}
\end{equation}
where $Z_0$ is the largest invariant subset of the set $Z = \{x\in D_0\,:\,\nabla  P(x)=0\}$ for the dynamical system~\eqref{discrete_DS}.
For an arbitrary $t\ge0$, denote  $N=\big[\frac{t}{\varepsilon}\big]$  and notice that  $0\le t-N{\varepsilon}<\varepsilon$.
Applying the triangle inequality together with Lemma~\ref{lemma_x}, we obtain:
$$
\begin{aligned}
  &\rho(x(t),S_0)=\inf_{y\in S_0}\|x(t)-y\|\\
  &\le\inf_{y\in S_0}\|x(N\varepsilon)-y\|+\|x(t)- x(N\varepsilon)\| \\
 &\le \inf_{y\in S_0}\|x(N\varepsilon)-y\|+\frac{M}{L}(e^{LU(N\varepsilon)\varepsilon}{-}1),
\end{aligned}
$$
where $U(N\varepsilon)=\max\limits_{s\in[N\varepsilon,t]}\sum_{i=1}^{\ell}|u^{\varepsilon}_{i}(s,N\varepsilon)|$.
From~\eqref{limit_discrete}, $\inf_{y\in S_0}\|x(N\varepsilon)-y\|\to 0$ as $N\to\infty$. Furthermore, \eqref{est_u}~yields
\begin{align*}
U(N\varepsilon)\varepsilon&\le C_1\varepsilon\|\nabla P(x(N\varepsilon))\|+{C_2}\big(\varepsilon\|\nabla P(x(N\varepsilon))\|\big)^{1/2}\\
&+{C_3}\big(\varepsilon\|\nabla P(x(N\varepsilon))\|\big)^{1/3},
\end{align*}
therefore, $U(N\varepsilon)\varepsilon\to 0$ as $N\to\infty$ because of the continuity of $\|\nabla P(x)\|$.
Summarizing all the above, we conclude that
$$
\rho(x(t),Z_0)\to 0 \text{ as }t\to\infty.
$$
%
%

%

\begin{thebibliography}{10}

\bibitem{Bro83}
{\sc Brockett, R.~W.}
\newblock Asymptotic stability and feedback stabilization.
\newblock {\em Differential Geometric Control Theory\/} (1983), 181--191.

\bibitem{Clar97}
{\sc Clarke, F.~H., Ledyaev, Y.~S., Sontag, E.~D., and Subbotin, A.~I.}
\newblock Asymptotic controllability implies feedback stabilization.
\newblock {\em IEEE Tran on Automatic Control 42}, 10 (1997), 1394--1407.

\bibitem{Jean14}
{\sc Jean, F.}
\newblock {\em Control of nonholonomic systems: from sub-{R}iemannian geometry
  to motion planning}.
\newblock Springer, 2014.

\bibitem{Kha86}
{\sc Khatib, O.}
\newblock Real-time obstacle avoidance for manipulators and mobile robots.
\newblock {\em The international journal of robotics research 5}, 1 (1986),
  90--98.

\bibitem{Kod90}
{\sc Koditschek, D.~E., and Rimon, E.}
\newblock Robot navigation functions on manifolds with boundary.
\newblock {\em Advances in Applied Mathematics 11}, 4 (1990), 412--442.

\bibitem{Kolman95}
{\sc Kolmanovsky, I., and McClamroch, N.~H.}
\newblock Developments in nonholonomic control problems.
\newblock {\em IEEE control systems 15}, 6 (1995), 20--36.

\bibitem{K09}
{\sc Kovalev, A., Martynyuk, A., Boichuk, O., Mazko, A., Petryshyn, R., Slyusarchuk, V., Zuyev, A., and Slyn’ko, V.}
\newblock Novel qualitative methods of nonlinear mechanics and their application to the analysis of multifrequency oscillations, stability, and control problems.
\newblock {\em Nonlinear Dynamics and Systems Theory 9}, 2 (2009), 117--145.

\bibitem{LaSalle}
{\sc LaSalle, J.}
\newblock {\em The stability and control of discrete processes}.
\newblock Springer, 1986.

\bibitem{LC90}
{\sc Li, Z., and Canny, J.}
\newblock Motion of two rigid bodies with rolling constraint.
\newblock {\em IEEE Transactions on Robotics and Automation 6}, 1 (1990),
  62--72.

\bibitem{Liu97}
{\sc Liu, W.}
\newblock An approximation algorithm for nonholonomic systems.
\newblock {\em SIAM Journal on Control and Optimization 35}, 4 (1997),
  1328--1365.

\bibitem{Pa16a}
{\sc Paternain, S., Koditschek, D.~E., and Ribeiro, A.}
\newblock Navigation functions for convex potentials in a space with convex
  obstacles.
\newblock {\em arXiv preprint arXiv:1605.00638\/} (2016).

\bibitem{Popa}
{\sc Popa, D.~O., and Wen, J.~T.}
\newblock Nonholonomic path-planning with obstacle avoidance: a path-space
  approach.
\newblock {\em Proc. 1996 IEEE Int. Conf. on Robotics and Automation 3\/}
  (1996), 2662--2667.

\bibitem{Rim92}
{\sc Rimon, E., and Koditschek, D.~E.}
\newblock Exact robot navigation using artificial potential functions.
\newblock {\em IEEE Tran on Robotics and Automation 8}, 5 (1992), 501--518.

\bibitem{Rou08}
{\sc Roussos, G.~P., Dimarogonas, D.~V., and Kyriakopoulos, K.~J.}
\newblock 3d navigation and collision avoidance for a non-holonomic vehicle.
\newblock {\em Proc. American Control Conf\/} (2008), 3512--3517.

\bibitem{Va12}
{\sc Sharma, B., Vanualailai, J., and Singh, S.}
\newblock Lyapunov-based nonlinear controllers for obstacle avoidance with a
  planar n-link doubly nonholonomic manipulator.
\newblock {\em Robotics and Autonomous Systems 60}, 12 (2012), 1484--1497.

\bibitem{Sic08}
{\sc Siciliano, B., and Khatib, O.}
\newblock {\em Springer Handbook of Robotics}.
\newblock Springer, 2008.

\bibitem{Tan01}
{\sc Tanner, H.~G., Loizou, S., and Kyriakopoulos, K.~J.}
\newblock Nonholonomic stabilization with collision avoidance for mobile
  robots.
\newblock {\em Proc. 2001 IEEE/RSJ Int Conf on Intelligent Robots and Systems
  3\/} (2001), 1220--1225.

\bibitem{Ura17}
{\sc Urakubo, T.}
\newblock Stability analysis and control of nonholonomic systems with potential
  fields.
\newblock {\em Journal of Intelligent \& Robotic Systems\/} (2017), 1--17, doi:
  10.1007/s10846--017--0473--1.

\bibitem{Va08}
{\sc Vanualailai, J., Sharma, B., and Nakagiri, S.}
\newblock An asymptotically stable collision-avoidance system.
\newblock {\em Int. J. Non-Linear Mechanics 43}, 9 (2008), 925--932.

\bibitem{Wal15}
{\sc Walters, P., Kamalapurkar, R., and Dixon, W.~E.}
\newblock Approximate optimal online continuous-time path-planner with static
  obstacle avoidance.
\newblock In {\em Proc. 54th IEEE Conf. on Decision and Control\/} (15-18
  December 2015), pp.~650--655.

\bibitem{YKD}
{\sc Yang, R., Krishnaprasad, P.~S., and Dayawansa, W.}
\newblock Optimal control of a rigid body with two oscillators.
\newblock In {\em Mechanics Day, Fields Institute Communications}. AMS, 1996,
  pp.~233--260.

\bibitem{Zu2001}
{\sc Zuyev, A.}
\newblock Application of control Lyapunov functions technique for partial stabilization.
\newblock In {\em IEEE Conference on Control Applications\/} (2001),
  pp.~509--513.

\bibitem{ZuSIAM}
{\sc Zuyev, A.}
\newblock Exponential stabilization of nonholonomic systems by means of
  oscillating controls.
\newblock {\em SIAM J Control Optim 54}, 3 (2016), 1678--1696.

\bibitem{ZuGrIJC}
{\sc Zuyev, A., and Grushkovskaya, V.}
\newblock Motion planning for control-affine systems satisfying low-order
  controllability conditions.
\newblock {\em International Journal of Control 90\/} (2017), 2517--2537.

\bibitem{ZuGrIFAC}
{\sc Zuyev, A., and Grushkovskaya, V.}
\newblock Obstacle avoidance problem for driftless nonlinear systems with
  oscillating controls.
\newblock {\em IFAC-PapersOnLine 50\/} (2017), 15343--15348.

\bibitem{ZuGrBeECC}
{\sc Zuyev, A., Grushkovskaya, V., and Benner, P.}
\newblock Time-varying stabilization of a class of driftless systems satisfying
  second-order controllability conditions.
\newblock In {\em Proc. 15th European Control Conference\/} (2016),
  pp.~1678--1696.

\end{thebibliography}
\end{document}